\documentclass[11pt]{article}
\usepackage{epsfig,amssymb}

\begin{document}

\title{
Gronwall's conjecture for 3-webs with  infinitesimal symmetries}

\author{{\Large Sergey I. Agafonov}\\
Departamento de Matem\'atica,\\
UNESP-Universidade Estadual Paulista,\\ S\~ao Jos\'e do Rio Preto, Brazil\\
e-mail: {\tt agafonov@ibilce.unesp.br} }
\date{}
\maketitle
\unitlength=1mm

\newtheorem{theorem}{Theorem}
\newtheorem{proposition}{Proposition}
\newtheorem{lemma}{Lemma}
\newtheorem{corollary}{Corollary}
\newtheorem{definition}{Definition}
\newtheorem{example}{Example}

\pagestyle{plain}

\begin{abstract}
We study non-flat planar 3-webs with infinitesimal symmetries. Using
multi-dimensional Schwarzian derivative we  give a criterion for linearization
of such webs and present a projective  classification thereof.
Using this classification we show that the Gronwall conjecture is true for 3-webs admitting infinitesimal symmetries.\\

\end{abstract}

\section{Introduction}
A planar 3-web $\mathcal{W}_3$  is formed by three foliations in the plane.
A 3-web  $\mathcal{L}_3$ is linear if the leaves of all its foliations are rectilinear.
 A 3-web $\mathcal{W}_3$ is
locally linearizable if there is a local diffeomorphism mapping $\mathcal{W}_3$ into some linear web $\mathcal{L}_3$. If this diffeomorphism can be chosen so
that the lines of each foliation of $\mathcal{L}_3$ are parallel, then the web is called
flat. Any projective transformation  maps a linear 3-web into a linear 3-web therefore a linearization, if it exists, is not unique.

Graf and Sauer (see \cite{GSg}) provided the following very elegant complete description of linear flat 3-webs.
Locally, to each of three foliations there corresponds a curve
 arc in the dual plane, thus we have 3 arcs. A linear web $\mathcal{L}_3$ is flat if and only if these 3 arcs belong to  some (possibly singular) cubic.
 Since any two flat 3-webs are locally diffeomorphic by definition, but not any two cubics are projectively equivalent, there are diffeomorphic but projectively
  non-equivalent linear 3-webs. For example, the 3-web formed by tangent lines to the curve dual to a smooth cubic  is not projectively
  equivalent to the web formed by three families of parallel lines,  whose corresponding cubic is decomposed in three lines.

 More than hundred years ago Gronwall conjectured in \cite{Gn} that any two local linearizations $\varphi,\psi$ of a non-flat 3-web
  are projectively equivalent, i.e. there exists a projective transformation  $g\in PGL(3,\mathbb C)$ such that $\psi=g\circ \varphi$. In spite of many efforts
  (see, for example, \cite{BB,Be,Bk,Bg,Bp,GLb,GMS,Hl,Su,Sa,Vn,Wg}) the conjecture is still open.
  Note that for a $d$-web with $d\ge 4$, there exists at most one projective equivalence class of linearizations,
  (see \cite{BB}).
  The reason is that a $4$-web determines a unique
projective connection such that the web leaves are geodesic and the linearizability of this web amounts to flatness of the connection
(see \cite{Lg,BB}) for the classical approach and \cite{Hl,AGL,Pl} for a modern exposition).

The treatment of the Gronwall conjecture in the literature is a little bit controversial.
Bol proved in \cite{Bk} that there are at most 17 essentially different (i.e. projectively non-equivalent)
 linearizations of a non-flat 3-web. Bor\r{u}vka in \cite{Bp} improved this bound to 16. In the short note \cite{Vn},
 Vaona  gave a sketchy proof that the bound is 11.
 Grifone, Muzsnay and Saab  studied the linearizability of 3-webs in \cite{GMS} and  proved that the bound is 15. Later
Goldberg and Lychagin (see \cite{GLb}) confirmed the result of \cite{GMS}, but claimed that the treatment of linearizability in \cite{GMS} is incomplete and the main example of
linearizable 3-web in \cite{GMS} is not linearizable.
Goldberg and Lychagin were also interested mainly in establishing linearizability criteria and
obtained their bound as a by-product.

There are also results on the Gronwall conjecture with some restrictions on the web and on the map.
For instance, Bol \cite{Bg} demonstrated the following two Theorems:
1) a local diffeomorphism mapping a linear 3-web, whose two families of  lines are tangents of some conic, to some 3-web of the same type is projective,
2) a local diffeomorphism mapping a linear 3-web, whose one family of  lines is a pencil  to some web of the same type, is projective.
Wang in \cite{Wg} proved that a 3-web admitting two projectively non-equivalent linearizations is flat,
provided that the Chern curvature of the web vanishes to order three at some point.
The reviewers of Mathematical Reviews and Zentralblatt MATH claim that Smirnov proved the Gronwall conjecture in a short note \cite{Su} published in an
obscure regional journal, which is virtually unavailable even in Russian scientific libraries. Note that in the same year Smirnov also presented
a more detailed paper \cite{Sa} in a respectable journal, where he proved Theorem 2 of Bol mentioned above.
Taking into account that Smirnov published several papers on the subject and understood the importance of the Gronwall conjecture,
we doubt that he would have published such a strong (and difficult) result under the title "On certain problems of uniqueness in the theory of webs"
in the journal of a provincial pedagogical university. We believe that there must be some misunderstanding on the part of the reviewers.
Actually, we will need only the fact that there are finitely many distinct linearizations.

This paper is motivated by the following result of Cartan \cite{Cg}: the symmetry algebra of a 3-web is either three-dimensional,
and then the web is flat; or one-dimensional; or trivial.
Thus the uniqueness of linearization does not hold true for the 3-webs admitting symmetry algebra of largest possible dimensions
and it is quite natural to look into the class of 3-webs with one-dimensional symmetry algebra.
(A similar strategy was successfully  applied in \cite{MPPi} to studying planar webs of maximal rank.)

To distinguish between essentially different linearizations one needs projective differential invariants of maps, namely
multi-dimensional Schwarzian derivatives (see \cite{Sp}).
Unlike the one-dimensional case, a complete invariant has more components than the map, thus leading to some differential relations between
the Schwarzian components.
These relations, which one can consider as differential syzygies, were used by H\'enaut
(see \cite{Hl}) to characterize the linearizability of planar webs in terms of solutions to these syzygies.
In fact, these equations, written for a special normalization of the vector fields tangent to the web leaves, were also the main tool of  Vaona in \cite {Vn}.
The principal difficulty in analyzing the linearizability of 3-webs lies in the fact that a 3-web does not determine uniquely a projective connection,
like in the case of 4-webs.
Therefore the "differential syzygies", manifesting the flatness of a "rectifying" projective connection,
give a weakly under-determined system of nonlinear PDEs. Compatibility conditions for this system are so hopelessly involved
 that are intractable even with the help of modern symbolic computation software.

If a linearizable non-flat 3-web has an infinitesimal symmetry, then
the Schwarzian derivative of the linearization is also invariant with respect to this symmetry. Thus we reduce the difficult problem of studying a weakly over-determined system of
non-linear PDEs to a simpler problem of compatibility of non-linear ODEs.
This compatibility condition amounts to vanishing of a resultant
of two polynomials, one being of degree 5 and the other of degree 6.

We also give a geometric characterization of one-dimensional symmetry of a non-flat 3-web.

Since there is a finite number of essentially different linearizations of a non-flat 3-web,
any infinitesimal symmetry of such web is necessarily projective.
This observation reduces the classification problem to the Jordan classification of $3\times 3$ matrices.

Using the obtained classification, we prove the Gronwall conjecture for linear 3-webs with infinitesimal symmetries. The proof is based on the fact
 that the conjecture is true for planar $d$-webs with $d>3$.
The key observations are the following:
 \begin{itemize}
  \item any of the classified 3-webs $\cal W$ admits an analytic extension to a global (singular) linear $d$-web $\widetilde{{\cal W}}$ on the projective plane
 with $d\ge 4$, the case $d=\infty$  (of countably many 1-parameter families of lines) being possible,
   \item if a diffeomorphism preserves linearity  of such 3-web $\cal W$ then it preserves linearity of some 4-subweb of the $d$-web $\widetilde{{\cal W}}$.
   \end{itemize}
This paper treats local properties of 3-webs, all the objects (webs, maps, infinitesimal symmetries) are defined in some open set $\Omega$.

\section{Non-projective morphisms of linear webs}\label{nonProj}
\begin{definition}
A morphism of a linear $d$-web is a local diffeomorphism that maps this $d$-web to a linear $d$-web.
\end{definition}
Let us recall some basic facts about projective differential invariants (for the detail see, for example, \cite{Sp}).
Two local diffeomorphisms $\varphi,\widetilde{\varphi}: \Omega\subset \mathbb C^2 \to \mathbb C^2$ are projectively equivalent,
i.e. there exists a projective transformation  $g\in PGL(3,\mathbb C)$ such that $\widetilde{\varphi}=g\circ \varphi$ if and only if their Schwarzian derivatives
$$S^k_{ij}=\widetilde{S^k_{ij}}$$
coincide, where  the components of the Schwarzian of the map $(u^1,u^2)\mapsto (\varphi ^1,\varphi^2)$ are defined as follows:
\begin{equation}\label{Schwarzian}
 S^k_{ij}(u)=\sum^2_{l=1}\frac{\partial ^2\varphi^l}{\partial u^i \partial u^j}\frac{\partial u^k}{\partial \varphi^l}
 -\frac{1}{3}\delta ^k_i \frac{\partial}{\partial u^j}\ln\left(\det\frac{\partial \varphi}{\partial u}\right)-\frac{1}{3}\delta ^k_j \frac{\partial}{\partial u^i}\ln\left(\det\frac{\partial \varphi}{\partial u}\right).
\end{equation}
They satisfy the following linear relations
$$
S^k_{ij}(u)=S^k_{ji}(u), \ \ \ \ \sum^2_{l=1}S^l_{il}(u)=0.
$$
To avoid working with a lot of indices let us choose the following notation for  linearly independent components of Schwarzian derivative:
\begin{equation}\label{KLMN}
K= S^1_{11},\ \ \ L= S^1_{22}, \ \ \ M= S^2_{11}, \ \ \ N= S^2_{22},
\end{equation}
and $x=u^1$, $y=u^2$ for the local coordinates.

Four functions $K,L,M,N$ are the above components of the Schwarzian derivative of some  map $\varphi: \Omega\subset \mathbb C^2 \to \mathbb C^2$ if and only if they satisfy the following non-linear PDEs:

\begin{equation}
\begin{array}{c}\label{compatibility}
2K_{xy}+M_{yy}+N_{xx}-6KK_y+2ML_x+LM_x+3NM_y-3KN_x+3MN_y=0,\\
 \\
K_{yy}+L_{xx}+2N_{xy}+3LK_x-3NK_y+3KL_x+ML_y+2LM_y-6NN_x=0.
\end{array}
\end{equation}
As "differential syzygies" of projective invariants, these equations were written explicitly by Tresse in \cite{Ti}.
In a slightly different but related context they appeared a bit earlier  in the Liouville studies of  projective connections \cite{Lg}.

Let $\mathcal{W}_3$ be a 3-web in some domain $\Omega\subset \mathbb C^2$ and $V_i= \partial_x+\lambda _i\partial_y$, $i=1,2,3$ vector fields tangent to the web leaves.
 The web is linear if and only if
  \begin{equation}\label{trans}
 V_i(\lambda_i)=0, \ \ \ i=1,2,3.
 \end{equation}
 H\'enaut in \cite{Hl} proved that $\mathcal{W}_3$ is linearizable by a diffeomorphism $\varphi=(\varphi^1,\varphi ^2)$ (i.e. $\varphi$ maps
  $\mathcal{W}_3$ into some linear 3-web $\mathcal{L}_3$) if and only if the components of the Schwarzian derivative of $\varphi$ verify
  \begin{equation}\label{cubH}
L\lambda^3_i-3N\lambda^2_i+3K\lambda_i-M= V_i(\lambda_i), \ \ \ i=1,2,3.
 \end{equation}
Since $\varphi$ is projective if and only if $K=L=M=N=0$, we have the following proposition.
\begin{proposition}
A linear 3-web admits a non-projective morphism if and only if there exists a non-vanishing solution to system (\ref{compatibility})
such that for the vector fields $V_i= \partial_x+\lambda _i\partial_y$,
$i=1,2,3$, tangent to the web leaves, holds true $L\lambda_i^3-3N\lambda_i^2+3K\lambda_i-M=0$.
\end{proposition}

\smallskip

\noindent {\bf Remark.} One can choose the local coordinates to normalize the vector fields to
$\partial _x,$ $\partial_y,$ $\partial _x+\lambda \partial_y.$ Then equations (\ref{cubH}) amount to $L=M=0$ and
$3\lambda(K-\lambda N)=\lambda_x+\lambda \lambda_y$ (for the detail see \cite{Hl}).

\section{Linearization of 3-webs with  infinitesimal symmetries}\label{linearization}
\begin{definition}
An infinitesimal symmetry of a d-web is a vector field whose local flow preserves the web.
 \end{definition}
All the infinitesimal symmetries of a web form a Lie algebra with respect to the Lie bracket.
As an example of his general theory of continuous transformations,  Cartan considered infinitesimal symmetries of 3-webs, formed by coordinate lines and  integral curves of one ODE in the plane,  and gave
 criteria for existing of nontrivial symmetry algebra. (See the original paper \cite{Cg} of Cartan. For a modern treatment the reader can look  up in \cite{GLb}.)
Following Cartan, we choose local coordinates $x,y$  so that the symmetry takes the form
\begin{equation}\label{NormalSymmetry}
Y=\partial _x+\partial_y,
\end{equation}
and the web leaves are tangent to the following vector fields
\begin{equation}\label{Vfields}
v_1=\partial _x, \ \ \ v_2=\partial_y, \ \ \ v_3=\partial _x+S(t)\partial_y,
\end{equation}
where $t=x-y$ is an invariant of the group action.
Then equations (\ref{cubH}) give $M=L=0$ and
\begin{equation}\label{Ocub}
3S(K-SN)=(1-S)S',
\end{equation}
where $S'=\frac{dS}{dt}$ (see the Remark in the end of Section \ref{nonProj}). Since there is only a finite number of projectively non-equivalent morphisms of a non-flat linear 3-web, one concludes immediately that $K,N$ are invariant with respect to the symmetry,
i.e. they are functions of $t$.  Now equations (\ref{compatibility})
take the form
\begin{equation}
\begin{array}{c}\label{Ocompatibility}
-2K''+N''+6KK'-3KN'=0,\\
 \\
K''-2N''+3NK'-6NN'=0.
\end{array}
\end{equation}
This system has the following two integrals:
\begin{equation}
\begin{array}{l}\label{Integrals}
I_1 =N'-K'+ K^2+N^2-NK,\\
 \\
I_2 = 3(N'K-K'N)-2(K^3+N^3)+3(K^2N+N^2K).
\end{array}
\end{equation}
Using the main result of \cite{Hl}, we conclude that the web is linearizable if and only if there is at least one solution to equations (\ref{Ocompatibility})
subject to relation (\ref{Ocub}). Due to the nonlinearity, the direct approach to finding the compatibility condition for the innocent looking overdetermined 
system (\ref{Ocub},\ref{Ocompatibility})
leads to very involved expressions. 
We simplify them by symmetrizing the system using the natural action of the symmetric group $S_3$, 
which permutes the vector fields (\ref{Vfields}). Transposing the vector fields $v_1$ and $v_2$ induces the transposition of $x$ and $y$ and the following action on our normalization of the web:
$$
\begin{array}{l}
x\mapsto y,\ \ y\mapsto x,\ \ t\mapsto -t,\\
S\mapsto \frac{1}{S},\ \ S'\mapsto \frac{S'}{S^2},\\
K\mapsto N,\ \ N\mapsto K.\\
\end{array}
$$
Similarly, transposing $y$, which is an integral of $v_1$, with $z$, the integral of $v_3$, defined by $dz=\frac{Sdx}{S-1}-\frac{dy}{S-1}$ to satisfy $Y(z)=1$, one gets, using  the formula for the Schwarzian derivative of a composition (see \cite{Sp}), the following substitutions:
$$
\begin{array}{l}
x\mapsto x,\ \ y\mapsto z,\ \ dt\mapsto \frac{dt}{1-S},\\
S\mapsto \frac{S}{S-1},\ \ S'\mapsto \frac{S'}{S-1},\\
K\mapsto K-2NS+\frac{2}{3}S',\ \ N\mapsto -N(S-1)+\frac{1}{3}S'.\\
\end{array}
$$
The above two transpositions generate the whole group action of $S_3$ that  permutes the fist integrals $x,y$ and $z$ of the foliations.

For a non-flat web holds true $\frac{d^2}{dt^2}(\log S)\ne 0$, therefore one can take the following symmetrization of $S$ as an independent variable:
\begin{equation}
X = \frac{1}{3}\frac{(S^2-S+1)^3}{(S-1)^2S^2},
\end{equation}
which is nothing else as the symmetrization of the cross-ratio of the directions of $Y,v_1,v_2,v_3$.
Symmetrizing $K$ and substituting $S'$ from equation (\ref{Ocub}) we define
\begin{equation}
U: = -\frac{K}{S-1}+\frac{SN}{S-1},
\end{equation}
similarly, applying the same procedure to $\frac{S^2K}{S+1}$ we define another invariant
\begin{equation}
V:=\frac{1}{3}\frac{(S^2-S+1)^2(N-K)}{(S-1)(S-2)(2S-1)(S+1)}.
\end{equation}
The first integrals (\ref{Integrals}) now take the forms
$$
\begin{array}{l}
\scriptstyle I_1 = 9U(4X-9)\frac{dV}{dX}+\frac{9(4X-9)}{X}V^2-\frac{6U(X-9)}{X}V+ U^2,\\
 \\
\scriptstyle I_2 =27U^2(4X-9)\frac{dV}{dX}+ \frac{27(4X-9)^2}{X^2}V^3-\frac{27U(4X-9)}{X}V^2-\frac{9U}{X}\left[2XU-18U+(12X^2-27X)\frac{dU}{dX}\right]V+2U^3.
\end{array}
$$
Note that the  the system (\ref{Ocub},\ref{Ocompatibility}) admits the symmetry algebra $\{\partial_t, -t\partial_t+K\partial_K+N\partial_N\}$. We have already used the first symmetry, choosing $X$ as a new variable. Now let us use the second one and set
$V=ZU$, $U'=FU$. Now $'$ means derivative by $X$. Substituting these expressions into the first integrals, taking total derivatives of them by $X$ and equating the results to zero, we get the following two equations:
\begin{equation}
\begin{array}{l}\label{eqZ}\scriptstyle
9X^2(4X-9)Z''+18X(4X-9)ZZ'+3X[9X(4X-9)F+10X+18]Z'  +[81+18X(4X-9)F]Z^2+\\
\\
\scriptstyle  [18X^2(4X-9)F^2+12X(9+2X)F-54+9X^2(4X-9)F']Z +  2X^2F=0,\\
  \\
\scriptstyle 9X^3(4X-9)Z''+27X(4X-9)^2Z^2Z'-18X^2(4X-9)ZZ'+3X^2[9X(4X-9)F+10X+18]Z'+\\
\\
\scriptstyle +27(4X-9)[6+X(4X-9)F]Z^3 -27X[3+X(4X-9)F]Z^2-18X[3+X(X-9)F]Z +2X^3F=0.\\

\end{array}
\end{equation}
For a non-flat web one has $F\ne 0$ since  due to (\ref{Ocub}) holds true $U = \frac{S'}{3S}$. Then a given non-flat 3-web has so many distinct linearizations as many solutions $Z(X)$ has the system (\ref{eqZ}).
Equations  (\ref{eqZ}) give $Z'$ and $Z''$ as functions of $Z,X,F,F'$. In particular,
 \begin{equation}\label{Zx}
\textstyle Z'= \frac{-9(4X-9)[6+X(4X-9)F]Z^2+3X[18+5X(4X-9)F]Z+X^2[6X(4X-9)F^2+(14X-18)F+3X(4X-9)F']}{3X(4X-9)[3(4X-9)Z-4X]}.
 \end{equation}
Expressions for $Z'$ and $Z''$ have the factor $3(4X-9)Z-4X$ in the denominators, but the equation $(12X-27)Z=4X$ is not compatible with (\ref{eqZ}).

The compatibility condition $\frac{dZ'}{dX}=Z''$ has the form $E(Z)=0$, where $E$ is a polynomal in $Z$ of degree 5
$$E(Z)=E_5Z^5+E_4Z^4+E_3Z^3+E_2Z^2+E_1Z+E_0,$$
 whose coefficients $E_i$ are given  in Appendix.

Equating the total derivative of $E(Z)$ to zero, one obtains a polynomial equation of degree 6 in $Z$
$$H(Z)=H_6Z^6+H_5Z^5+H_4Z^4+H_3Z^3+H_2Z^2+H_1Z+H_0=0,$$
where the coefficients $H_i$ are given in Appendix.  Let us define $R(X,F,F',F'',F''')$ as the resultant of polynomials $E(Z)$ and $H(Z)$, the functions $\rho,\ \omega$ by
$$\rho (X,F,F'):=X(4X-9)^2F'+2X(4X-9)^2F^2+6(X-1)(4X-9)F-8,$$ $$\omega(X,F,F',F'',F'''):=\frac{R(X,F,F',F'',F''')}{(4X-9)^{20}X^{26}\rho ^8(X,F,F')}.$$
\begin{theorem}
A non-flat 3-web admitting  infinitesimal symmetry is linearizable if and only if the invariant $F(X)$ satisfies $\omega(X,F,F',F'',F''')=0$.
\end{theorem}
{\it Proof:} Our 3-web is linearizable if and only if the equations $E(Z)=0$, $H(Z)=0$ have a common solution $Z(X)$,
i.e. if and only if the resultant of $E(Z)$ and $H(Z)$ vanishes. The resultant factors out as $(4X-9)^{20}X^{26}\rho ^8\omega$.
If $\rho (X,F,F')=0$ then one finds $F'$ from this equation  and computes $F''$ and $F'''$.
With this expressions for the derivatives of $F$, the polynomials  $E(Z)$ and $H(Z)$ have the common factor $\sigma (X,Z)=[3(4X-9)Z-4X]^3$.
As was mentioned above, the equation $\sigma (X,Z)=0$ is not compatible with  system  (\ref{eqZ}) for $Z$.
Now the resultant of $\frac{E(Z)}{\sigma (X,Z)}$ and $\frac{H(Z)}{\sigma (X,Z)}$ is the polynomial in $F$:
$$
\scriptstyle \tau(F):=36X^6(4X-9)^6[(84X-189)F-20][7X^2(4X-9)^3F^3+21X(51X-4)(4X-9)^2F^2-12(596X-21)(4X-9)F+11856].
$$
A direct computation shows that $\tau(F)=0$ is not compatible with $F'$ defined by $\rho (X,F,F')=0$.
\hfill $\Box$

\noindent{\bf Remark 1.} The computation in the proof of the above theorem was made with the help of symbolic computation software, namely Maple. The polynomial $\omega$ is quintic in the highest derivative $F'''$, the coefficient of $(F''')^5$ being equal to $X^{11}(4X-9)^8$. The expressions for the other coefficients are hopelessly involved to be presented here in their generality.\\
 \\
\noindent{\bf Remark 2.} One is tempted to apply the Euclid algorithm to $H(Z)$ and $E(Z)$ to prove that there is at most one linearization. It works only for the first step, i.e. one can check by Maple that there is at most 4 linearizations, but the next step is out of Maple's reach.\\
 \\
\noindent{\bf Remark 3.} Normalization (\ref{Vfields}) is defined by the web up to permutations
 of the first integrals $x,y,z$ of foliations and up to scaling and translating of the parameter  $t$.
These transformations generate a group $G$ acting on parameterized curves $(t,S(t))$, differential expressions $X,F$ being invariants of this action.
(Note that $U = \frac{S'}{3S}$ and therefore $F$ do not depend on the linearizing diffeomorphism.)
Then the parameterized curve $(X,F(X))$ is an analog of the signature curve considered by Olver (see \cite{Oi}) for the  action of the Euclidean group in the plane.\\

\noindent Finally, let us give a geometric characterization of one-dimensional symmetry of a non-flat 3-web (for a analytical criterion of existence in terms of differential invariants see, for instance, \cite{GLb}).
\begin{theorem}
A non-flat 3-web formed by integral curves of three vector fields $v_1,v_2,v_3$ has an infinitesimal symmetry $Y$ if and only if the three 3-webs generated similarly by the  3-tuples of vector fields $\{Y,v_2,v_3\}$, $\{v_1,Y,v_3\}$, and $\{v_1,v_2,Y\}$ are flat.
\end{theorem}
{\it Proof:} If a non-flat 3-web has a symmetry $Y$ then, using the Cartan normalization (\ref{NormalSymmetry},\ref{Vfields}), one easily sees that the 3-web generated by $\{v_1,v_2,Y\}$ is flat. This proves the flatness of the three 3-webs of  $\{Y,v_2,v_3\}$, $\{v_1,Y,v_3\}$, and $\{v_1,v_2,Y\}$.

To prove the converse claim note that the direction of the vector field $Y$ cannot be tangent to the web leaves; otherwise one normalizes the vector fields to take the form $Y=v_1=\partial _x$, $v_2=\mu(y)\partial _x+\partial _y$, $v_3=\nu(y)\partial _x+\partial _y$ and our 3-web is obviously flat. Let us again choose the coordinates so that  the vector field $Y$ takes the form  (\ref{NormalSymmetry}) and $v_1,v_2$ are as in (\ref{Vfields}). This is possible due to the flatness of the web $\{v_1,v_2,Y\}$.   Now the 3-web $\{v_1,v_2,v_3=\partial _x+\lambda \partial_y\}$ can be defined by three 1-forms as follows

$$
\omega _1=dy,\ \ \ \ \omega _2=dx,\ \ \ \ \omega _3=dy-\lambda dx.
$$
Further, the flatness of the web $\{v_1,Y,v_3\}$ implies
$$
\lambda(\lambda-1)(\lambda_{xx}+\lambda_{xy})=(2\lambda-1)\lambda_x(\lambda_x+\lambda_y).
$$
Similarly, the flatness of $\{Y,v_2,v_3\}$ gives
$$
(\lambda-1)(\lambda_{xy}+\lambda_{yy})=\lambda_y(\lambda_x+\lambda_y).
$$
Computing the compatibility conditions for the above two equations for $\lambda$, we arrive at
$$
(\lambda_x+\lambda_y)[\lambda_x^2(1-2\lambda)-\lambda \lambda_x\lambda_y+\lambda_{xx}\lambda(\lambda-1)]=0.
$$
If the factor in the square brackets vanishes then one finds all the second order derivatives of $\lambda$, 
in particular $\lambda_{xy}=\frac{\lambda_x\lambda_y}{\lambda}$, which implies that the web of $\{v_1,v_2,v_3\}$ is flat. 
Thus this factor is not zero and therefore $\lambda_x+\lambda_y=0$. That means that the foliation by the integral curves of the equation $\omega_3=0$ is also invariant by $Y$.
\hfill $\Box$

\section{Linear 3-webs with one-dimensional infinitesimal symmetry}\label{classes}
In this section we give a classification of non-flat 3-webs with infinitesimal symmetries up to projective transformations.
A symmetry algebra is called projective if it generates a subgroup of the projective group.

\begin{lemma}\label{isprojective}
If a linear non-flat 3-web admits a one-dimensional infinitesimal symmetry algebra then this symmetry algebra  is projective.
\end{lemma}
{\it Proof:} Let  $g^t$ be the local flow of the symmetry. For each $t$ the map $g^t$ is a morphism of our web. Since there is only a finite number of projectively non-equivalent morphisms for a non-flat linear 3-web (see \cite{Bk}), one has  $g^t \in PGL(3,\mathbb C)$.
\hfill $\Box$\\
\smallskip

Let $\mathcal{L}_3$ be a linear 3-web with a one-dimensional infinitesimal symmetry. Its $i$th foliation is a one-parameter family of lines in the form
\begin{equation}\label{DFC}
y=p_i(t)x+q_i(t),
\end{equation}
or the pencil $x=cst$.
 Thus we obtain parameterized curves $(p_i(t),q_i(t))$  or a segment of the line at infinity $l_{\infty}$ in the dual plane.
 In what follows we call them dual focal curves and denote them by $\Phi_i$.
 Recall that a 3-web $\mathcal{L}_3$ is flat if and only if these three dual focal curves are arcs of a same cubic curve (Graf and Sauer Theorem \cite{GSg}).
 Our projective symmetry acts also  in the dual space. Obviously, the dual focal curves are (locally) invariant.

In the dual plane, the action of the projective group is generated by the following operators (see, for instance, \cite{Oe}):
\begin{equation}
\begin{array}{l}
T_1=\partial_p,\ \ \ T_2=\partial_q,\ \ \ A_1=p\partial_q,\ \ \ A_2=q\partial_p,\ \ \ D_1=p\partial_p,\\
\\
D_2=q\partial_q,\ \ \ \Pi_1=p^2\partial_p+pq\partial_q,\ \ \ \ \Pi_2=pq\partial_p+q^2\partial_q,
\end{array}
\end{equation}
where $p,q$ are affine coordinates of the line $y=px+q$.

Now let us classify one-dimensional subalgebras of the Lie algebra $pgl(3,\mathbb C)$.
This classification is provided by the Jordan normal forms of $3\times 3$ matrices with zero traces.
One can also normalize one of the non-zero eigenvalues to 1.
Thus one can take the following matrices as the orbit representatives:
$$\scriptstyle
\Xi_1=\left(\begin{array}{ccc}\scriptstyle
 0 &\scriptstyle 1 & \scriptstyle 0 \\
 \scriptstyle 0 & \scriptstyle 0 & \scriptstyle 1\\
 \scriptstyle 0 & \scriptstyle 0 & \scriptstyle 0 \\
 \end{array}\right), \ \
 \Xi_{2,1}=\left(\begin{array}{ccc}
 \scriptstyle 0 & \scriptstyle 0 & \scriptstyle 0 \\
 \scriptstyle 0 &\scriptstyle  0 &\scriptstyle  1\\
 \scriptstyle 0 & \scriptstyle 0 &\scriptstyle 0 \\
 \end{array}\right),  \ \
 \Xi_{2,3}=\left(\begin{array}{ccc}
 \scriptstyle 1 &\scriptstyle 1 &\scriptstyle 0 \\
 \scriptstyle 0 &\scriptstyle 1 &\scriptstyle 0\\
 \scriptstyle 0 &\scriptstyle 0 &\scriptstyle -2 \\
 \end{array}\right),
 \ \
 \Xi_{3,2}=\left(\begin{array}{ccc}
 \scriptstyle 1 &\scriptstyle 0 &\scriptstyle 0 \\
 \scriptstyle 0 &\scriptstyle -1 &\scriptstyle 0\\
 \scriptstyle 0 &\scriptstyle 0 &\scriptstyle 0 \\
 \end{array}\right)
 \ \
 \Xi_{3,3}=\left(\begin{array}{ccc}
\scriptstyle a &\scriptstyle 0 &\scriptstyle 0 \\
\scriptstyle 0 &\scriptstyle b &\scriptstyle 0\\
\scriptstyle 0 &\scriptstyle 0 &\scriptstyle -(a+b) \\
 \end{array}\right),
$$
where $a \ne 0$, $b \ne 0$, $a+b \ne 0$.
The first subindex is the number of Jordan's blocks and the second (if any) is the matrix rank.
Let $l=(P,Q,R)$ be homogeneous coordinates of a line $RY=PX+QZ$. If the matrix representation of the infinitesimal action is $M l^T$, then in the affine coordinates the operator is
\begin{equation}\label{matrixtodiff}
\begin{array}{c}
  \{(M_{11}-M_{33})p+M_{12}q + M_{13}-M_{31}p^2-M_{32}pq\}\partial_p+ \\
 \{M_{21}p+(M_{22}-M_{33})q+M_{23}-M_{31}pq-M_{32}q^2\}\partial_q.\\
\end{array}
\end{equation}
For example, the operator corresponding to the matrix $\Xi_1$ is $q\partial_p+\partial_q$. We take as a representative of the same orbit the operator
$\xi _1=\partial_p+p\partial_q$. The invariant curves for this operator are the ones parameterized by $\left(p,p^2/2+\lambda\right)$ with some constant
$\lambda$ and the line at infinity $l_{\infty}$. To the line $l_{\infty}$ there corresponds the foliation by parallel lines $x=cst$.
We can move  these curves around by the stabilizer of the algebra $\{\xi_1\}$. This stabilizer is spanned by the following operators: $\{T_2,T_1+A_1,D_1+2D_2\}$.
As a line intersects an invariant curve at most at 2 points, the dual focal curves $\Phi_i$ either are pieces of pairwise
 distinct invariant curves or two $\Phi_i$'s are pieces of the same invariant curve.
 In the first case, using the stabilizer we can bring them to the following parameterized forms:
 $$(p_1,p_1^2/2),\ \ (p_2,p_2^2/2+1),\ \ (p_3,p_3^2/2+\lambda),\ \ \lambda \neq 0,1,$$
 and in the second case one sets $\lambda=0$.
 The corresponding linear 3-webs are given by the following families of lines:
$$y=p_1x+\frac{p_1^2}{2},\ \ \ y=p_2x+\frac{p_2^2}{2}+1,\ \ \  y=p_3x+\frac{p_3^2}{2}+\lambda.$$
Finally, if one of the dual focal curves is included in $l_{\infty}$, we can parameterized the other two either as follows:
$$(p_1,p_1^2/2),\ \ (p_2,p_2^2/2+1),$$
whose 3-web is
$$y=p_1x+\frac{p_1^2}{2},\ \ \ y=p_2x+\frac{p_2^2}{2}+1,\ \ \  x+p_3=0,$$
or
$$(p_1,p_1^2/2),\ \ (p_2,p_2^2/2),$$ with the hexagonal web
\begin{equation}\label{Sigma1_3}
y=p_1x+\frac{p_1^2}{2},\ \ \ y=p_2x+\frac{p_2^2}{2},\ \ \  x+p_3=0.
\end{equation}
To write the corresponding symmetry operator $\check{\xi}$ in "geometric" coordinates we note that the passage to the dual plane is a contact transform given by the following formulae:
$$
\begin{array}{c}
  p=\frac{dy}{dx},  \ \ \ q=y-\frac{dy}{dx}x, \ \ \ \frac{dq}{dp}=-x,\\
  x= -\frac{dq}{dp}, \ \ \ y=q-\frac{dq}{dp}p,  \ \ \ \frac{dy}{dx}=p.
\end{array}
$$
Therefore one has $\check{\xi} _1=-\partial_x+x\partial_y$

Proceeding with this scheme, one obtains the classification below, where the parameters $\lambda,\mu$ in the parametrization of dual focal curves are supposed to be pairwise distinct
and different from the corresponding constants in the normalized curves.
\begin{theorem}
Any non-flat linear 3-web with infinitesimal symmetry is projectively equivalent to one in Table \ref{class1}. Moreover, these normal forms are projectively nonequivalent.
\end{theorem}
\begin{table}
\noindent \begin{tabular}{|c|c|c|c|}
\hline

Type    & Operator                             &            Dual Focal Curves          \\
\hline
$\Xi_1$ & $\begin{array}{c}
             \partial_p+p\partial_q, \\
             \\
             -\partial_x+x\partial_y
             \end{array}$
                                               &

                                                                         $\begin{array}{ll}
                                                                             &\\
                                                                          1) & (p_1,p_1^2/2)  \\
                                                                             & (p_2,p_2^2/2+1) \\
                                                                             & (p_3,p_3^2/2+\lambda)\\
                                                                             &\\
                                                                         \end{array} $
                                                                         $\begin{array}{ll}
                                                                         &\\
                                                                          2)&(p_1,p_1^2/2)  \\
                                                                             &(p_2,p_2^2/2) \\
                                                                             &(p_3,p_3^2/2+1)\\
                                                                             &\\
                                                                         \end{array} $
                                                                          $\begin{array}{ll}
                                                                         &\\
                                                                          3)&(p_1,p_1^2/2)  \\
                                                                             &(p_2,p_2^2/2+1) \\
                                                                             &l_{\infty}\\
                                                                             &\\
                                                                         \end{array} $                                                  \\

\hline

$\Xi_{3,2}$ & $\begin{array}{c}
                 p\partial_p-q\partial_q, \\
                 \\
                 -2x\partial_x-y\partial_y
               \end{array}$                        &

                                                                         $\begin{array}{ll}
                                                                            &\\
                                                                          1) & (p_1,1/p_1)  \\
                                                                             & (p_2,\lambda/p_2) \\
                                                                             & (p_3,\mu/p_3)\\
                                                                             &\\
                                                                          2)& (p_1,1/p_1)  \\
                                                                             & (p_2,1/p_2) \\
                                                                             & (p_3,\lambda/p_3)\\
                                                                             &\\
                                                                             \end{array} $
                                                                         $\begin{array}{ll}
                                                                             &\\
                                                                          3)& (p_1,1/p_1)  \\
                                                                             & (p_2,\lambda/p_2) \\
                                                                             & (p_3,0)\\
                                                                             &\\
                                                                          4)& (p_1,1/p_1)  \\
                                                                             & (p_2,0) \\
                                                                             & (0,1/p_3)\\
                                                                             &\\
                                                                         \end{array} $
                                                                         $\begin{array}{ll}
                                                                             &\\
                                                                          5)& (p_1,1/p_1)  \\
                                                                             & (p_2,\lambda/p_2) \\
                                                                             & l_{\infty}\\
                                                                             &\\
                                                                          6)& (p_1,1/p_1)  \\
                                                                             & (p_2,0) \\
                                                                             & l_{\infty}\\
                                                                             &\\
                                                                         \end{array} $
                                                                                                                            \\

\hline
$\Xi_{2,3}$ & $\begin{array}{c}
                 p\partial_p+(p+q)\partial_q, \\
                 \\
                 -\partial_x+y\partial_y
               \end{array}$                   &

                                                                         $\begin{array}{ll}
                                                                             &\\
                                                                          1) & (p_1,p_1\ln p_1)  \\
                                                                             & (p_2,p_2[\ln p_2+\lambda]) \\
                                                                             & (p_3,p_3[\ln p_3+\mu])\\
                                                                             &\\
                                                                          2)& (p_1,p_1\ln p_1)  \\
                                                                             & (p_2,p_2\ln p_2) \\
                                                                             & (p_3,p_3[\ln p_3+\lambda])\\
                                                                             &\\
                                                                          3)& (p_1,p_1\ln p_1)  \\
                                                                             & (p_2,p_2\ln p_2) \\
                                                                             & (p_3,p_3\ln p_3)\\
                                                                             &\\
                                                                         \end{array} $
                                                                         $\begin{array}{ll}
                                                                         &\\
                                                                          4)& (p_1,p_1\ln p_1)  \\
                                                                             & (p_2,p_2[\ln p_2+\lambda]) \\
                                                                             & (0,p_3)\\
                                                                             &\\
                                                                           5)& (p_1,p_1\ln p_1)  \\
                                                                             & (p_2,p_2\ln p_2) \\
                                                                             & (0,p_3)\\
                                                                             & \\
                                                                          6) & (p_1,p_1\ln p_1)  \\
                                                                             & (p_2,p_2[\ln p_2+\lambda]) \\
                                                                             & l_{\infty}\\
                                                                             &\\
                                                                         \end{array} $
                                                                         $\begin{array}{ll}
                                                                         &\\
                                                                          7)& (p_1,p_1\ln p_1)  \\
                                                                             & (p_2,p_2\ln p_2) \\
                                                                             & l_{\infty}\\
                                                                             &\\
                                                                           8)& (p_1,p_1\ln p_1)  \\
                                                                             & (0,p_2) \\
                                                                             & l_{\infty}\\
                                                                             & \\
                                                                             & \\
                                                                             & \\
                                                                             & \\
                                                                             &\\
                                                                         \end{array} $
                                                                                                                            \\
\hline
$\Xi_{3,3}$ & $\begin{array}{c}
                 p\partial_p+\beta q\partial_q, \\
                 \\
                 (\beta-1)x\partial_x+y\partial_y\\
                 \\
                 \beta \ne 0,1,-1,2,\frac{1}{2}
               \end{array}$                              &
                                                                         $\begin{array}{ll}
                                                                             &\\
                                                                          1) & (p_1,p_1^\beta)  \\
                                                                             & (p_2,\lambda p_2^\beta) \\
                                                                             & (p_3,\mu p_3^\beta)\\
                                                                             &\\
                                                                          2)& (p_1,p_1^\beta)  \\
                                                                             & (p_2,p_2^\beta) \\
                                                                             & (p_3,\lambda p_3^\beta)\\
                                                                             &\\
                                                                          3)& (p_1,p_1^\beta)  \\
                                                                             & (p_2,p_2^\beta) \\
                                                                             & (p_3,p_3^\beta)\\
                                                                             & \beta \ne 3,\frac{1}{3}, \frac{3}{2}, \frac{2}{3}\\
                                                                             &\\
                                                                         \end{array} $
                                                                         $\begin{array}{ll}
                                                                             &\\
                                                                         4)& (p_1,p_1^\beta)   \\
                                                                             & (p_2,\lambda p_2^\beta) \\
                                                                             & (p_3,0)\\
                                                                             &\\
                                                                           5)& (p_1,p_1^\beta)  \\
                                                                             & (p_2,p_2^\beta) \\
                                                                             & (p_3,0)\\
                                                                             &\\
                                                                           6)& (p_1,p_1^\beta)  \\
                                                                             & (p_2,0) \\
                                                                             & (0,p_3)\\
                                                                             &\\
                                                                             &\\
                                                                         \end{array} $
                                                                         $\begin{array}{ll}
                                                                             &\\
                                                                          7)& (p_1,p_1^\beta)   \\
                                                                             & (p_2,\lambda p_2^\beta) \\
                                                                             & l_{\infty}\\
                                                                             &\\
                                                                           8)& (p_1,p_1^\beta)  \\
                                                                             & (p_2,p_2^\beta) \\
                                                                             & l_{\infty}\\
                                                                             &\\
                                                                          9)& (p_1,p_1^\beta)  \\
                                                                             & (p_2,0) \\
                                                                             & l_{\infty}\\
                                                                             &\\
                                                                             &\\
                                                                         \end{array} $
                                                                                                                            \\
\hline
\end{tabular}
\caption{Classification of non-flat symmetric 3-webs}\label{class1}
\end{table}

\smallskip

\noindent {\it Proof:}  To the matrix $\Xi_{2,1}$ there corresponds the operator $\partial_q$ whose stabilizer is spanned by $\{T_1,T_2,A_1,D_1,D_2\}$; the invariant curves are rectilinear
(lines $p=cst$ and $l_{\infty}$) and the web is flat.

The invariant curves of the operator $p\partial_p-q\partial_q$, corresponding to $\Xi_{3,2}$, are hyperbolas $pq=cst$ and the lines $p=0$,  $q=0$,
$l_{\infty}$. The stabilizer is spanned by $\{D_1,D_2\}$. Here the web is flat if and only if its two dual focal curves are arcs of the same hyperbola and
the third is a segment of one of the lines or if and only if all its dual focal curves are rectilinear. 
The forms in the table are clearly not projectively equivalent since the lines $p=0$,  $q=0$ are tangent to the hyperbolas, 
whereas the line at infinity $l_{\infty}$  cuts them at 2 points.

For the symmetry type $\Xi_{2,3}$ we choose the operator $p\partial_p+(p+q)\partial_q$ as a representative
(one easily gets this operator applying formula (\ref{matrixtodiff}), transposing $p,q$, and then scaling $q$).
The invariant curves are $q/p-\ln p=cst$ and the lines $p=0$, $l_{\infty}$.
Here $\ln$ is the multivalued analytical function, the above formula being the parametrization of the line family by points of the Riemann surface of $\ln$.
The stabilizer is spanned by $\{A_1,D_1+D_2\}$.
A projective transform, mapping the line at infinity $l_{\infty}$ in the line $p=0$ does not preserve the symmetry hence
the forms in the table are not projectively equivalent.

To the symmetry type $\Xi_{3,3}$ we get the operator $\alpha p\partial_p+\beta q\partial_q$, where $\alpha=2a+b$ and $\beta=a+2b$.
Therefore $[\alpha:\beta]\ne [1:2],[2:1],[1:-1]$.
As the cases $\alpha=0$ (or $\beta=0$) and $\alpha=\beta$ give  flat 3-webs, one can choose $\alpha=1$ and $\beta \ne 0,1,-1,2,\frac{1}{2}$.
The invariant curves are $q=\lambda p^{\beta}$ (here again $p^{\beta}$ is multivalued for non-integer $\beta$) or lines $p=0$, $q=0$, $l_{\infty}$.
The stabilizer is spanned by $\{D_1,D_2\}$.
The invariant curve is a  cubic  if and only if $\beta \in \{3,\frac{1}{3},\frac{3}{2},\frac{2}{3}\}$, therefore these values are excluded for the case 3,
where the dual focal curves are included in the same invariant curve. (Recall the Graf and Sauer Theorem.)
Again, a projective transform, mapping the line at infinity $l_{\infty}$ in the line $p=0$ (or $q=0$), does not preserve the symmetry hence
 the forms in the table are not projectively equivalent.
\hfill $\Box$\\

\noindent{\bf Remark.} As a by-product of the above proof,  we obtain normal forms up to projective equivalence for the triple of dual focal curves
of flat linear 3-webs with projective infinitesimal symmetries, as well as the corresponding  algebras of infinitesimal simmetries of this type.
Namely, one can choose them as follows:
\begin{enumerate}
  \item the cuspidal cubic $p=q^3$, one-dimensional algebra $\{3p\partial_p+q\partial_q\},$
  \item the conic $pq=1$ and its secant $l_{\infty}$, one-dimensional algebra $\{ p\partial_p-q\partial_q\},$
  \item the conic $2q=p^2$ and its tangent line $l_{\infty}$, two-dimensional algebra $\{\partial_p+p\partial_q, p\partial_p+2q\partial_q\},$
  \item three non-concurrent lines $p=0$, $q=0$, $l_{\infty}$,  two-dimensional algebra $\{ p\partial_p, q\partial_q\},$
  \item three concurrent lines $p=0$, $p=1$, $p=-1$,  three-dimensional algebra $\{\partial_q, p\partial_q, q\partial_q\}.$
\end{enumerate}

\section{Gronwall's conjecture for  3-webs with infinitesimal symmetries}
The following theorem implies that the Gronwall conjecture  is true for  3-webs with infinitesimal symmetries.
\begin{theorem}
The normal forms in Table \ref{class1} are pairwise not diffeomorphic. Moreover, any diffeomorphism, preserving a normal form, is projective.
\end{theorem}
{\it Proof:} Let us fix some normal form. The three dual focal curves $\Phi_i$ are pieces of three curves invariant under the symmetry.
Some of this invariant curves are permitted to coincide, but for each normal form a generic line in the dual space intersects this collection of invariant curves
in more then 3 points.  This is obvious for the symmetry types $\Xi_1$, $\Xi_{3,2}$ and $\Xi_{3,3}$ with a rational value of the parameter $\beta$,
where the invariant curves are algebraic.
For the the symmetry type $\Xi_{2,3}$, a generic line
$q=kp+l,\ \ k,l=cst$  intersects the invariant curve $q/p-\ln p=cst$ in infinitely many points: substituting $p=e^z$, one gets this conclusion from the theorem on
values of a holomorphic function in a neighborhood of essential singularity, the singularity being $z=\infty$.
For the type $\Xi_{3,3}$ with an irrational  $\beta$, one applies the same trick with the  substitution $p=e^z$.

Therefore any of the classified 3-webs $\cal W$ admits an analytic extension to a global (singular) linear $d$-web $\widetilde{{\cal W}}$ on the projective plane
 with $d\ge 4$, the case $d=\infty$  (of countably many 1-parameter families of lines) being possible.
The key observation: if a local diffeomorphism  preserves linearity of our 3-web then it also preserves linearity of some 4-subweb of that $d$-web.
Hence the diffeomorphism is a projective transform, since the Gronwall conjecture  is true for 4-webs.

Below we present the scheme for proving the observation and  work out all the details for the web in the normal form $\Xi_1,3$,
which will be denoted by ${\cal W}(\Xi_1,3)$.

The leaves of each foliation of the web are locally parameterized. At a non-singular point,
where the leaves are transverse, one can choose two of these parameters as local coordinates.
Then the parameter of the third family is a function of the chosen ones,
determined explicitly by some equation relating all three parameters.
This relation is called the web equation of the web. (Note that its form depends on the choice of parameters.)
The web equation, relating the chosen parameters $t_i$ along the dual focal curves (\ref{DFC}),   reads as
\begin{equation}\label{web equation}
\det \left(
       \begin{array}{ccc}
         1 & p_1 & q_1 \\
         1 & p_2 & q_2 \\
         1 & p_3 & q_3 \\
       \end{array}
     \right)=0.
\end{equation}
The geometric meaning of this equation is that the three lines (from different foliations) corresponding to the parameters $t_i$, satisfying the equation, are concurrent.
If there are $k$ values, say, of $t_3$ satisfying this equation for fixed $t_1,t_2$ and giving different points $(p_3,q_3)$ in the dual plane by virtue of (\ref{DFC}),
 then there are $k$ lines of the third family passing through the intersection
point of   the two lines, one from the first family with $t_1$ and the other of the second one with $t_2$, and our local web can be extended to $d$-web with
$d\ge k+2$.

\noindent As the local parameters for ${\cal W}(\Xi_1,3)$, let as choose $p_1$ and $p_3$, then the web equation reads as
$$
p_2^2-2p_3p_2+(2+2p_1p_3-p_1^2)=0.
$$
The values of $p_3$ and $p_1$ fix the point $(x,y)$, one of the solutions $p_2$ to this quadratic web equation gives the third line of the web.
(Note that this equation defines, in fact, a 5-web $\widetilde{{\cal W}}(\Xi_1,3)$: there are 2 solutions for $p_2$ and there are 2 lines of the first family passing
through $(x,y)$.)

 Web equation locally defines some surface $M_1$ in 3-dimensional space with the parameters $t_1,t_2,t_3$ as local coordinates.
 Moving the lines of the web, the infinitesimal symmetry generates also an action on this surface $M_1$. Moreover, this action is the restriction of
 the local flow of some vector field in 3-dimensional space of parameters. For the web ${\cal W}(\Xi_1,3)$ this vector field is
 $\partial _{p_1}+\partial _{p_2}+\partial _{p_3}$, for the other normal forms the operators are presented in the column "Symmetry" of Table \ref{data}.

This vector field has two first integrals. Being invariant under the symmetry, the web equation can be written in terms of these two invariants.
This defines  a curve in a two-dimensional
space. Thus to each of the normal forms there corresponds a Riemann surface $S_1$ (one-dimensional analytic manifold). For the web ${\cal W}(\Xi_1,3)$,
choosing the  invariants $w=p_2-p_1, \ \ z= p_3-p_1$ one gets   the Riemann surface $S_1$:
$$
w^2-2zw+2=0
$$
as the symmetry reduction  of the web equation. The invariants for the other normal forms one finds in the column "Invariants" of Table \ref{data}.
The equations, defining the corresponding Riemann surfaces are given in the column "Riemann surfaces" of Table \ref{data}.

\noindent Note that the parameters $z,w$ on these Riemann surfaces are chosen so that:
 \begin{enumerate}
 \item each pair of lines $l_1,l_2$ of our linear 3-web fixes a value of $z$ and, by duality, two points on the invariant curves,
 \item $S_1$ becomes the   Riemann surface  of a multi-valued analytic function $\widetilde{w}(z)$, defined by some equation $f(z,w)=0$
 (see the column "Riemann surfaces" of Table \ref{data} for the explicit formulas for $f$),
 \item one of the values of $\widetilde{w}(z)$ gives the third line $l_3$, passing through the intersection point of $l_1$ and $l_2$,
 \item this third line $l_3$ defines the third point on one of the invariant curves.

 \end{enumerate}
 Suppose a local diffeomorphism maps the chosen linear 3-web from the
 list  to some (possibly the same) linear 3-web from the list. This fixes a pair of symmetry types $(\Xi_{*},\Xi_*)$,
 where the first element is the symmetry type of the first web and the second is the symmetry type of the second web.

This diffeomorphism induces a map between the corresponding surfaces $M_1$, $M_2$,
defined by the web equations of the webs. On the surfaces $M_1$ and $M_2$ our webs are represented as the 3-webs cut by the planes $t_i=cst$, $T_j=cst,$
where $T_j$ are parameters along the dual focal curves of the second web.
Mapping  the 3-web on $M_1$ to the 3-web on $M_2$, the diffeomorphism takes the form $T_{\pi(i)}=g_i(t_i)$,
where  $\pi$ is some permutation of 3 indices.
Moreover, the diffeomorphism relates the corresponding symmetry operators.   This condition gives ODEs of the first order for $g_i$.
Resolving these ODEs and taking into account that the symmetry operator on the first surface $M_1$
 is mapped to some multiple of the symmetry operator on the second surface $M_2$, one finds $g_i$  up to 4 constants.

Consider a diffeomorphism of the web ${\cal W}(\Xi_1,3)$, say, to the web ${\cal W}(\Xi_1,1)$ whose web equation  and symmetry are
$$
(P_2-P_1)P_3^2+(P_1^2-P_2^2-2)P_3+[P_1P_2^2-P_1^2P_2+2\lambda (P_2-P_1)+2P_1]=0
$$
and $\partial _{P_1}+\partial _{P_2}+\partial _{P_3}$ respectively.
Permuting the indices, if necessary,  we conclude that the diffeomorphism takes the form $P_1=g_1(p_1), \ \ P_2=g_3(p_3)$.
Since the symmetry is preserved, we have
$$
g_1(p_1)=kp_1+c_1, \ \ g_3(p_3)=kp_3+c_3,
$$
 where $k,c_1,c_3$ are some constants. Moreover, under this diffeomorphism holds true $P_3=kp_2+c_2$ for some constant $c_2$.
 For the other pairs of symmetry types, the formulas for $g_i$ are presented below.

Further, our diffeomorphism between the surfaces $M_1$ and $M_2$ maps the orbits of the symmetry on $M_1$
 to the orbits of the symmetry on $M_2$.
Therefore the diffeomorphism is lowerable to a local biholomorphism of the corresponding Riemann surfaces
$S_1$ defined by $f(z,w)=0$ and $S_2$ defined by $F(Z,W)=0$. In particular, this biholomorphism has the form $Z=a(z),\ W=b(w)$. Using the explicit formulas for the functions $g_i$,
one easily checks that $a(z),b(w)$ either are analytic for any $z,w\in \mathbb C$ or have a branch point at zero.
Thus $F(a(z),b(w(z)))\equiv 0$ on some neighborhood of $(z_0,w_0)\in S_1$.

The corresponding symmetry reduction   of the web equation for ${\cal W}(\Xi_1,1)$ is the Riemann surface $S_2$:
$$
ZW^2-(Z^2+2)W+2\lambda Z=0,
$$
where $Z=P_2-P_1, \ \ W=P_3-P_1.$
Our local diffeomorphism induces the following local biholomorphism
\begin{equation}\label{bih}
Z=kz+(c_3-c_1), \ \ W=kw+(c_2-c_1).
\end{equation}

Using the formulas  from the column "Riemann surfaces" of Table \ref{data}, one  easily checks that
the locally defined analytic function $w(z)$ is extended to some multi-valued function $\tilde{w}(z)$.

\begin{lemma}\label{path}
For each of the following pairs of symmetry types $(\Xi_1,\Xi_1)$,
 $(\Xi_{3,2},\Xi_1)$, $(\Xi_{3,2},\Xi_{3,2})$,
 $(\Xi_{2,3},\Xi_1)$, $(\Xi_{2,3},\Xi_{3,2})$, $(\Xi_{2,3},\Xi_{2,3})$,
 $(\Xi_{3,3},\Xi_1)$, $(\Xi_{3,3},\Xi_{3,2})$, $(\Xi_{3,3},\Xi_{2,3})$,$(\Xi_{3,3},\Xi_{3,3})$,
one can choose a closed path $\gamma \subset \mathbb C, \ \ z_0\in \gamma$  so that:
\begin{enumerate}
\item  the analytic continuation of $a(z)$ along $\gamma$ does not change the branch,
\item  the analytic continuation of $w(z)$ along $\gamma$ changes the branch of $\tilde{w}(z)$,
\item  the value $w_1$ of $\tilde{w}(z)$ on this new branch defines a 4th line passing through the intersection point of the concurrent lines $l_1,l_2,l_3$.
\end{enumerate}
\end{lemma}
The proof the Lemma is presented after the proof of the Theorem.

Since $F(a(z),b(w(z)))\equiv 0$ on some neighborhood of $z_0$ and $F(a(z),b(w(z)))$ is analytic,
we conclude that $F(a(z),b(\tilde{w}(z)))\equiv 0$ along the path.
Thus the same map $Z=a(z),\ W=b(w)$ maps some neighborhood $\Sigma_1$ of $(z_0,w_1)\in S_1$, where $w_1$ is the value of $\tilde{w}(z)$ on the new branch,
to same neighborhood $\Sigma_2$ of the Riemann surface $S_2$. The neighborhood $\Sigma_1$ defines locally some "additional" family of lines,
one of which passes through the point determined by $z_0$, and similarly, $\Sigma_2$ defines locally some family of lines,
one of which passes through the point determined by $Z_0$.
This means that the diffeomorphism maps the "additional" family of lines of the first $d$-web $\widetilde{{\cal W}}$
to some "additional" family of lines of the second one.

For the web ${\cal W}(\Xi_1,3)$, let us choose a closed path $\gamma$ with $z_0 \in \gamma $ so that it goes around one of the branch points $z_b=\pm \sqrt 2$ of $S_1$.
(Here $z_0$ corresponds to the base point of our web.) When we come back to $z_0$, we change the branch of $\tilde{w}(z)$.
Along the path holds true $F(kz+(c_3-c_1),kw+(c_2-c_1))\equiv 0$, where $F(Z,W)=ZW^2-(Z^2+2)W+2\lambda Z=0$ defines the Riemann surface $S_2$.
Thus to the value $w_1=\overline{p_2}-p_1$ on the new branch corresponds a line different from $l_1,l_2,l_3$,
since $w_1=\overline{p_2}-p_1\ne w=p_2-p_1$ implies $\overline{p_2}\ne p_2$.
Therefore some 4th family of lines of the web $\widetilde{{\cal W}}(\Xi_1,3)$ is mapped by the diffeomorphism to some family of lines
of the web $\widetilde{{\cal W}}(\Xi_1,1)$, which is in fact 6-web.
(Note that the dual focal curves of the web $\widetilde{{\cal W}}(\Xi_1,1)$ are three different parabolas
 and a generic line intersects them in 6 points.)
Since any morphism of linear 4-web is projective, our diffeomorphism (if it exists) is also projective. But the forms ${\cal W}(\Xi_1,1)$ and ${\cal W}(\Xi_1,3)$ are not projectively equivalent.

Now let us prove that any diffeomorphism preserving ${\cal W}(\Xi_1,3)$ is projective. Again this diffeomorphism is of the form
$$
P_{\pi(i)}=kp_i+c_i, \ \ i=1,2,3.
$$
It generates an automorphism (\ref{bih}) of $S_1$, where the new parameters on $S_1$ are chosen as  $W=P_{\pi(2)}-P_{\pi(1)}, \ \ Z= P_{\pi(3)}-P_{\pi(1)}$. Repeating the trick with analytic continuation along the closed path we conclude that the diffeomorphism is projective.
\hfill $\Box$\\

\bigskip

\noindent {\it Proof of Lemma \ref{path}:} The case of the pair $(\Xi_1,\Xi_1)$ was considered in the proof of the Theorem: the function  
$a(z)$ is one-valued,
the path $\gamma $ goes around one of the branch points of $\widetilde{w}(z)$, which are not zero.

\begin{table}
\noindent \begin{tabular}{|c|c|c|c|l|}
\hline
Type      & Symmetry                & Invariants         & Riemann Surface                     & Branch Point \\
\hline
$\Xi_1,1$  & $\sum_i\partial _{p_i}$ & $\begin{array}{c}
                                         z=p_2-p_1 \\
                                         w=p_3-p_1  \\
                                         \end{array}$      &  $zw^2-(z^2+2)w+2\lambda z=0$       &  $z_b^4+(4-8\lambda)z_b^2+4=0$   \\
\hline
$\Xi_1,2$  & $\sum_i\partial _{p_i}$ & $\begin{array}{c}
                                        z=p_2-p_1 \\
                                        w=p_3-p_1  \\
                                        \end{array}$      &  $w^2-zw+2=0$       &  $z_b^2-8=0$   \\
\hline
$\Xi_1,3$  & $\sum_i\partial _{p_i}$ & $\begin{array}{c}
                                        z=p_3-p_1 \\
                                        w=p_2-p_1  \\
                                        \end{array}$      &  $w^2-2zw+2=0$       &  $z_b^2-2=0$   \\
\hline
$\Xi_{3,2},1$  & $\sum_ip_i\partial _{p_i}$ & $\begin{array}{c}
                                               z=p_2/p_1\\
                                               w=p_3/p_1  \\
                                               \end{array}$      &   $\begin{array}{l}
                                                                      (z-\lambda)w^2+(\lambda-z^2)w+\\
                                                                      -\mu z(z-1)=0
                                                                      \end{array}$                        &
                                                                                                                     $\begin{array}{l}
                                                                                                                     z_b^4+4\mu z_b^3+\\
                                                                                                                     -2(\lambda+2\mu \lambda+2\mu)z_b^2+\\
                                                                                                                     4\mu \lambda z_b+\lambda^2=0  \\
                                                                                                                      \end{array}$ \\
\hline
$\Xi_{3,2},2$  & $\sum_ip_i\partial _{p_i}$ & $\begin{array}{c}
                                               z=p_2/p_1\\
                                               w=p_3/p_1  \\
                                               \end{array}$      &   $w^2-(z+1)w-\lambda z=0$ & $\begin{array}{l}
                                                                                                                    z_b^2+(4\lambda+2)z_b+1=0\\
                                                                                                                      \end{array}$ \\
\hline
$\Xi_{3,2},3$  & $\sum_ip_i\partial _{p_i}$ & $\begin{array}{c}
                                               z=p_3/p_1\\
                                               w=p_2/p_1  \\
                                               \end{array}$      &   $w^2-wz+\lambda (z-1)=0$ & $\begin{array}{l}
                                                                                                                    z_b^2-4\lambda z_b+4\lambda=0\\
                                                                                                                      \end{array}$ \\
\hline
$\Xi_{3,2},4$  & $\sum_ip_i\partial _{p_i}$ & $\begin{array}{c}
                                               z=p_3/p_2\\
                                               w=p_1/p_2  \\
                                               \end{array}$      &   $w^2-w+z=0$ & $\begin{array}{l}
                                                                                                                    1-4z_b=0\\
                                                                                                                      \end{array}$ \\
\hline
$\Xi_{3,2},5$  & $\begin{array}{l}
                \sum_i^2p_i\partial _{p_i}+\\
                 -2p_3\partial _p3
                 \end{array}$
                                              & $\begin{array}{c}
                                               z=p_3p_2^2\\
                                               w=p_1/p_2  \\
                                               \end{array}$      &   $zw^2+(\lambda-z)w-1=0$ & $\begin{array}{l}
                                                                                                                    z_b^2+(4-2\lambda)z_b+\lambda ^2=0\\
                                                                                                                      \end{array}$ \\
\hline
$\Xi_{3,2},6$  & $\begin{array}{l}
                \sum_i^2p_i\partial _{p_i}+\\
                 -2p_3\partial _p3
                 \end{array}$
                                              & $\begin{array}{c}
                                               z=p_3p_2^2\\
                                               w=p_1/p_2  \\
                                               \end{array}$      &   $zw^2-zw-1=0$ & $\begin{array}{l}
                                                                                                                    z_b+4=0\\
                                                                                                                      \end{array}$ \\
\hline
$\Xi_{2,3},1$  & $\begin{array}{l}
                \sum_i\partial _{t_i},\\
                 p_i=e^{t_i}
                 \end{array}$
                                              & $\begin{array}{c}
                                               z=t_1-t_3\\
                                               w=t_2-t_3  \\
                                               \end{array}$      &   $\begin{array}{l}
                                                                       we^w(1-e^z)+ze^z(e^w-1)+ \\
                                                                       (\lambda-\mu) e^w+\mu e^z-\lambda e^{z+w}=0
                                                                       \end{array}$                                 & $\begin{array}{l}
                                                                                                                    \\
                                                                                                                      \end{array}$ \\
\hline
$\Xi_{2,3},2$  & $\begin{array}{l}
                \sum_i\partial _{t_i},\\
                 p_i=e^{t_i}
                 \end{array}$
                                              & $\begin{array}{c}
                                               z=t_1-t_3\\
                                               w=t_2-t_3  \\
                                               \end{array}$      &   $\begin{array}{l}
                                                                       we^w(1-e^z)+ze^z(e^w-1)+ \\
                                                                       \lambda (e^z- e^w)=0
                                                                       \end{array}$                                 & $\begin{array}{l}
                                                                                                                    \\
                                                                                                                      \end{array}$ \\
\hline
$\Xi_{2,3},3$  & $\begin{array}{l}
                \sum_i\partial _{t_i},\\
                 p_i=e^{t_i}
                 \end{array}$
                                              & $\begin{array}{c}
                                               z=t_1-t_3\\
                                               w=t_2-t_3  \\
                                               \end{array}$      &   $\begin{array}{l}
                                                                       we^w(1-e^z)+ze^z(e^w-1)=0\\
                                                                       \end{array}$                                 & $\begin{array}{l}
                                                                                                                    \\
                                                                                                                      \end{array}$ \\
\hline
$\Xi_{2,3},4$  & $\begin{array}{l}
                \sum_i\partial _{t_i},\\
                 p_i=e^{t_i}
                 \end{array}$
                                              & $\begin{array}{c}
                                               z=t_1-t_3\\
                                               w=t_2-t_3  \\
                                               \end{array}$      &   $\begin{array}{l}
                                                                       e^w-e^z+(w-z+\lambda)e^{z+w}=0\\
                                                                       \end{array}$                                 & $\begin{array}{l}
                                                                                                                    \\
                                                                                                                      \end{array}$ \\
\hline
$\Xi_{2,3},5$  & $\begin{array}{l}
                \sum_i\partial _{t_i},\\
                 p_i=e^{t_i}
                 \end{array}$
                                              & $\begin{array}{c}
                                               z=t_1-t_3\\
                                               w=t_2-t_3  \\
                                               \end{array}$      &   $\begin{array}{l}
                                                                       e^w-e^z+(w-z)e^{z+w}=0\\
                                                                       \end{array}$                                 & $\begin{array}{l}
                                                                                                                    \\
                                                                                                                      \end{array}$ \\
\hline
$\Xi_{2,3},6$  & $\begin{array}{l}
                \sum_i\partial _{t_i}+\partial _{p_3},\\
                 p_i=e^{t_i}, \\
                  i=1,2
                 \end{array}$
                                              & $\begin{array}{c}
                                               z=t_1-p_3\\
                                               w=t_2-p_3  \\
                                               \end{array}$      &   $\begin{array}{l}
                                                                      we^w-ze^z+\lambda e^{w}=0\\
                                                                       \end{array}$                                 & $\begin{array}{l}
                                                                                                                    \\
                                                                                                                      \end{array}$ \\
\hline
$\Xi_{2,3},7$  & $\begin{array}{l}
                \sum_i\partial _{t_i}+\partial _{p_3},\\
                 p_i=e^{t_i}, \\
                  i=1,2
                 \end{array}$
                                              & $\begin{array}{c}
                                               z=t_1-p_3\\
                                               w=t_2-p_3  \\
                                               \end{array}$      &   $\begin{array}{l}
                                                                      we^w-ze^z=0\\
                                                                       \end{array}$                                 & $\begin{array}{l}
                                                                                                                    \\
                                                                                                                      \end{array}$ \\
\hline
$\Xi_{2,3},8$  & $\begin{array}{l}
                \sum_i\partial _{t_i}+\partial _{p_3},\\
                 p_i=e^{t_i}, \\
                  i=1,2
                 \end{array}$
                                              & $\begin{array}{c}
                                               z=t_2-p_3\\
                                               w=t_1-p_3  \\
                                               \end{array}$      &   $\begin{array}{l}
                                                                      we^w-e^z=0\\
                                                                       \end{array}$                                 & $\begin{array}{l}
                                                                                                                    \\
                                                                                                                      \end{array}$ \\
\hline
$\Xi_{3,3},1$  & $\begin{array}{l}
                \sum_i\partial _{t_i},\\
                 p_i=e^{t_i}
                 \end{array}$
                                              & $\begin{array}{c}
                                               z=t_2-t_1\\
                                               w=t_3-t_1  \\
                                               \end{array}$      &   $\begin{array}{l}
                                                                       \mu e^{{\beta}w}(e^z-1)-\lambda e^{{\beta}z}(e^w-1)+\\
                                                                       e^w-e^z=0
                                                                       \end{array}$                                 & $\begin{array}{l}
                                                                                                                    \\
                                                                                                                      \end{array}$ \\
\hline
$\Xi_{3,3},2$  & $\begin{array}{l}
                \sum_i\partial _{t_i},\\
                 p_i=e^{t_i}
                 \end{array}$
                                              & $\begin{array}{c}
                                               z=t_2-t_1\\
                                               w=t_3-t_1  \\
                                               \end{array}$      &   $\begin{array}{l}
                                                                       \lambda e^{{\beta}w}(e^z-1)- e^{{\beta}z}(e^w-1)+\\
                                                                       e^w-e^z=0
                                                                       \end{array}$                                 & $\begin{array}{l}
                                                                                                                    \\
                                                                                                                      \end{array}$ \\
\hline
$\Xi_{3,3},3$  & $\begin{array}{l}
                \sum_i\partial _{t_i},\\
                 p_i=e^{t_i}
                 \end{array}$
                                              & $\begin{array}{c}
                                               z=t_2-t_1\\
                                               w=t_3-t_1  \\
                                               \end{array}$      &   $\begin{array}{l}
                                                                       e^{{\beta}w}(e^z-1)- e^{{\beta}z}(e^w-1)+\\
                                                                       e^w-e^z=0
                                                                       \end{array}$                                 & $\begin{array}{l}
                                                                                                                    \\
                                                                                                                      \end{array}$ \\
\hline
\end{tabular}
\caption{Reduction of web equations.}\label{data}
\end{table}
\setcounter{table}{1}
\begin{table}
\noindent \begin{tabular}{|c|c|c|c|l|}
\hline
Type      & Symmetry                & Invariants         & Riemann Surface                     & Branch Point \\
\hline
$\Xi_{3,3},4$  & $\begin{array}{l}
                \sum_i\partial _{t_i},\\
                 p_i=e^{t_i}
                 \end{array}$
                                              & $\begin{array}{c}
                                               z=t_1-t_3\\
                                               w=t_2-t_3  \\
                                               \end{array}$      &   $\begin{array}{l}
                                                                       \lambda e^{{\beta}w}(e^z-1)- e^{{\beta}z}(e^w-1)=0\\

                                                                       \end{array}$                                 & $\begin{array}{l}
                                                                                                                    \\
                                                                                                                      \end{array}$ \\
\hline
$\Xi_{3,3},5$  & $\begin{array}{l}
                \sum_i\partial _{t_i},\\
                 p_i=e^{t_i}
                 \end{array}$
                                              & $\begin{array}{c}
                                               z=t_1-t_3\\
                                               w=t_2-t_3  \\
                                               \end{array}$      &   $\begin{array}{l}
                                                                       e^{{\beta}w}(e^z-1)- e^{{\beta}z}(e^w-1)=0\\

                                                                       \end{array}$                                 & $\begin{array}{l}
                                                                                                                    \\
                                                                                                                      \end{array}$ \\
\hline
$\Xi_{3,3},6$  & $\begin{array}{l}
                \sum_i\partial _{t_i},\\
                 p_i=e^{t_i}\\
                 p_3=e^{{\beta}t_3}
                 \end{array}$
                                              & $\begin{array}{c}
                                               z=t_2-t_3\\
                                               w=t_1-t_3  \\
                                               \end{array}$      &   $\begin{array}{l}
                                                                       e^z(e^{{\beta}w}-1)+e^w=0\\

                                                                       \end{array}$                                 & $\begin{array}{l}
                                                                                                                    \\
                                                                                                                      \end{array}$ \\
\hline
$\Xi_{3,3},7$  & $\begin{array}{l}
                \sum_i\partial _{t_i},\\
                 p_i=e^{t_i}\\
                 p_3=e^{({\beta}-1)t_3}
                 \end{array}$
                                              & $\begin{array}{c}
                                               z=t_1-t_3\\
                                               w=t_2-t_3  \\
                                               \end{array}$      &   $\begin{array}{l}
                                                                  \lambda e^{{\beta}w}-e^w+e^z-e^{{\beta}z}=0\\
                                                                       \end{array}$                                 & $\begin{array}{l}
                                                                                                                    \\
                                                                                                                      \end{array}$ \\
\hline
$\Xi_{3,3},8$  & $\begin{array}{l}
                \sum_i\partial _{t_i},\\
                 p_i=e^{t_i}\\
                 p_3=e^{({\beta}-1)t_3}
                 \end{array}$
                                              & $\begin{array}{c}
                                               z=t_1-t_3\\
                                               w=t_2-t_3  \\
                                               \end{array}$      &   $\begin{array}{l}
                                                                  e^{{\beta}w}-e^w+e^z-e^{{\beta}z}=0\\
                                                                       \end{array}$                                 & $\begin{array}{l}
                                                                                                                    \\
                                                                                                                      \end{array}$ \\
\hline
$\Xi_{3,3},9$  & $\begin{array}{l}
                \sum_i\partial _{t_i},\\
                 p_i=e^{t_i}\\
                 p_3=e^{({\beta}-1)t_3}
                 \end{array}$
                                              & $\begin{array}{c}
                                               z=t_1-t_2\\
                                               w=t_3-t_2  \\
                                               \end{array}$      &   $\begin{array}{l}
                                                                  e^{({\beta}-1)w}(1-e^z)+e^{{\beta}z}=0\\
                                                                       \end{array}$                                 & $\begin{array}{l}
                                                                                                                    \\
                                                                                                                      \end{array}$ \\
\hline
\end{tabular}
\caption{Reduction of web equations.}\label{data2}
\end{table}

A diffeomorphism of a web with the symmetry type $\Xi_{3,2}$ to some web with the symmetry type  $(\Xi_1,1)$ would have the form
$$
P_{\pi(i)}=k\ln(p_i)+c_i, \ \ i=1,2,3,
$$
except for the forms 5 and 6, where $P_{\pi(3)}=-\frac{k}{2}\ln(p_3)+c_3$.
For such diffeomorphisms we choose a path $\gamma$ going around  one of the branch points of $\widetilde{w}(z)$ but not around the origin.

A diffeomorphism of a web with the symmetry type $\Xi_{3,2}$ to some web with the same symmetry would have the form
$$
P_{\pi(i)}=c_ip_i^k, \ \ i=1,2,3.
$$
(For the forms 5 and 6 one adjusts the exponent in the obvious way.)
Again we choose a path $\gamma$ going around a branch point of $\widetilde{w}(z)$ but not around the origin.

In the above considered cases of the pairs of symmetry types, the claim of the third item of Lemma \ref{path} is obviously true, since the invariant $w$ is linear
in the inclination  $p$ of web lines.

A diffeomorphism of a web with the symmetry type $\Xi_{2,3}$ to some web with the symmetry type $\Xi_1$ would have the form
$$
P_{\pi(i)}=kt_i+c_i, \ \ i=1,2,3.
$$
(For the forms 6,7,8 one substitutes $t_3=p_3$.) Here $a(z)$ is one-valued,
we choose a path $\gamma $ going around one of the branch points of $\widetilde{w}(z)$.

A diffeomorphism of a web with the symmetry type $\Xi_{2,3}$ to some web with the symmetry type $\Xi_{3,2}$ would be (up to substitution $t_3=p_3$ for the forms 6,7,8)
$$
P_{\pi(i)}=c_ie^{kt_i}, \ \ i=1,2,3.
$$
(For the forms $\Xi_{3,2},5$ and  $\Xi_{3,2},6$  one again adjusts the exponent.) Again $a(z)$ is one-valued,
the path $\gamma $ goes around one of the branch points of $\widetilde{w}(z)$.

A diffeomorphism of a web with the symmetry type $\Xi_{2,3}$ to some web with the same symmetry type  would be (up to substitution $t_3=p_3$ for the forms 6,7,8)
$$
T_{\pi(i)}=kt_i+c_i, \ \ i=1,2,3.
$$
Again $a(z)$ is one-valued,
the path $\gamma $ goes around one of the branch points of $\widetilde{w}(z)$.

To prove the third item of Lemma \ref{path} note that: 1) parameter $t_i$ defining the line $l_4$ is related to the inclination $p$  of web lines by $p_i=e^{t_i}$,
2) the value $w_1$ on the new branch cannot be related to the value $w$ defining the line $l_3$ by the formula $w_1=w+2\pi n$, $n\in \mathbb Z$,
since from the equations defining the Riemann surface $S_1$ one easily sees that  $(z,w)\in S_1\  \Rightarrow (z,w+2\pi n)\notin S_1 $.

Now we have all the formulas written and the reader can easily check the statement of the Lemma  for the pairs $(\Xi_{3,3},\Xi_{*})$.
For irrational $\beta$ one proves the third item of Lemma \ref{path} exactly as in the case of pairs $(\Xi_{2,3},\Xi_{*})$.
For rational $\beta$ one observes that $(\overline{Z},\overline{W})=(e^z,e^{\widetilde{w}(z)})$ satisfies some algebraic equation and defines an algebraic function
 $\overline{W}=\mathcal{A}(\overline{Z})$. Now one chooses the path $\gamma$ so that $\mathcal{A}(e^z)$ changes the branch, which ensures that $l_4$ is different from
  $l_1,l_2,l_3$.
We present equations for the branch point $z_b$ only for the forms with the symmetry type $\Xi_1$ and for the cases when it is important that $z_b\ne 0$.
\hfill $\Box$

\section{Acknowledgement}
The author thanks D.V.~Alekseevsky for useful discussions.
This research was partially supported by FAPESP grant \#2014/17812-0.
The author is also very grateful to the reviewers for their valuable comments and suggestions, which helped to improve the manuscript.

\section{Appendix: coefficients of $E(Z)$ and $H(Z)$}
$$
\begin{array}{l}\label{coef5}
\scriptstyle E_5=81(4X-9)^3,\\
\\
\scriptstyle E_4=-324X(4X-9)^2, \\
\\
\scriptstyle E_3=-27X^2(4X-9)[X(4X-9)^2F'+2X(4X-9)^2F^2+ (6X-13)(4X-9)F-12],\\
\\
\scriptstyle E_2=-3X^2(4X-9)[3X^2(4X-9)^2F''+21X^2(4X-9)^2FF'+6X(4X-9)(2X-9)F'+\\
 \\
 \scriptstyle 18X^2(4X-9)^2F^3+ 42X(4X-9)(X-3)F^2+ (-108X^2+144X+54)F-8], \\
  \\
\scriptstyle E_1=3X^3[8X^2(4X-9)^2F''+ 53X^2(4X-9)^2FF'+ 2X(40X-69)(4X-9)F'+42X^2(4X-9)^2F^3+\\
\\
\scriptstyle 2X(95X-153)(4X-9)F^2+ (108-180X)F-16],\\
 \\
\scriptstyle E_0=X^4[-16X^2(4X-9)F'' +3X^2(4X-9)^2(F')^2+12X^2(4X-9)^2F^2F'- 4X(16X+9)(4X-9)FF'- \\
\\
\scriptstyle  24X(8X-9)F'+ 12X^2(4X-9)^2F^4-72X(4X-9)F^3+ (108+240X-300X^2)F^2-(48X+16)F],\\
\\
\scriptstyle
H_6=-243(4X-9)^4[5X(4X-9)F-12X+30],\\
\\
\scriptstyle
H_5=81X(4X-9)^3[73X(4X-9)F-192X+486],\\
\\
\scriptstyle  H_4= 27X^2(4X-9)^2[-3X^2(4X-9)^3F''-3X^2(4X-9)^3FF'-9X(10X-21)(4X-9)^2F'+18X^2(4X-9)^3F^3+\\
\\
\scriptstyle -3X(30X-61)(4X-9)^2F^2-2(4X-9)(180X^2-578X+747)F+1008X-2592], \\
  \\
\scriptstyle H_3= 9X^2(4X-9)[-3X^3(4X-9)^4F'''+ -15X^3(4X-9)^4FF''-6X^2(14X-33)(4X-9)^3F''-\\
\\
\scriptstyle 21X^3(4X-9)^4(F')^2 -12X^3(4X-9)^4F^2F' -3X^2(215X-384)(4X-9)^3FF'+ 6X(30X^2+121X-360)(4X-9)^2F'+\\
 \\
\scriptstyle 36X^3(4X-9)^4F^4-102X^2(5X-6)(4X-9)^3F^3 -9X(102X^2-461X+534)(4X-9)^2F^2 +\\
\\
\scriptstyle 2(4X-9)(1800X^3-6548X^2+5445X+810)F+2160+2712X-1344X^2],  \\
\\
\end{array}
$$

$$
\begin{array}{l}\label{coef6}

\scriptstyle H_2=-3(4X-9)X^3[-36X^3(4X-9)^3F'''-189X^3(4X-9)^3FF''-18X^2(72X-133)(4X-9)^2F''-\\
 \\
\scriptstyle 216X^3(4X-9)^3(F')^2 -117X^3(4X-9)^3F^2F'-18X^2(424X-737)(4X-9)^2FF' -\\
\\
\scriptstyle 72X(4X-9)(90X^2-370X+357)F'+
414X^3(4X-9)^3F^4 -18X^2(221X-359)(4X-9)^2F^3+\\
\\
\scriptstyle  -18X(4X-9)(1002X^2-3484X+2985)F^2 + (8640X^3-22040X^2-5328X+16524)F-2448+960X], \\
\\
\scriptstyle H_1= -X^4[144X^3(4X-9)^3F'''+792X^3(4X-9)^3FF''+ 48X^2(116X-201)(4X-9)^2F''+ \\
\\
\scriptstyle 54X^3(4X-9)^4F(F')^2+18x

^2(34X+15)(4X-9)^3(F')^2+216X^3(4X-9)^4F^3F'+ 9X^2(61X+48)(4X-9)^3F^2F'+\\
\\
\scriptstyle 12X(2336X^2-3921X-270)(4X-9)^2FF'+36X(4X-9)(1040X^2-3560X+2709)F'+ 216X^3(4X-9)^4F^5-\\
\\
\scriptstyle 18X^2(43X+12)(4X-9)^3F^4+ 18X(619X^2-1027X-252)(4X-9)^2F^3+ \\
\\
\scriptstyle 24(4X-9)(3054X^3-9961X^2+7146X+405)F^2 +(11520X^3-83264X^2+148536X-27864)F+ 6048-2304X], \\
\\
\scriptstyle H_0= X^5[64X^3(4X-9)^2F'''+368X^3(4X-9)^2FF''+160X^2(16X-27)(4X-9)F''+\\
\\
\scriptstyle 63X^3(4X-9)^3F(F')^2+ 14X^2(8X+27)(4X-9)^2(F')^2 +252X^3(4X-9)^3F^3F' +4X^2(64X+189)(4X-9)^2F^2F'+ \\
\\
\scriptstyle 8X(4X-9)(1336X^2-2025X-567)FF'+16X(300X-269)(4X-9)F'+252X^3(4X-9)^3F^5+192X^3(4X-9)^2F^4+\\
\\
\scriptstyle 4X(7X+9)(95X-189)(4X-9)F^3+(13608-86904X^2+41904X+28800X^3)F^2+(3840X^2-7840X-2016)F
].

\end{array}
$$

\end{document}